\documentclass[10pt, a4paper, leqno]{article}
\usepackage[dvips]{hyperref}
\usepackage[british]{babel}
\usepackage{enumerate}
\usepackage{amsmath, amscd, amsfonts, amssymb, latexsym, theorem}
\usepackage{xypic}
\usepackage[T1]{fontenc}
\usepackage[latin1]{inputenc}

{\theorembodyfont{\itshape}   \newtheorem{thm}{Theorem}[section]}
{\theorembodyfont{\itshape}   \newtheorem{lem}[thm]{Lemma}}
{\theorembodyfont{\itshape}   }
{\theorembodyfont{\itshape}   \newtheorem{rem}[thm]{Remark}}
{\theorembodyfont{\itshape}   \newtheorem{prop}[thm]{Proposition}}
{\theorembodyfont{\itshape}   \newtheorem{cor}[thm]{Corollary}}
{\theorembodyfont{\itshape}   \newtheorem{ass}[thm]{Assumption}}

\DeclareMathOperator{\Sym}{Sym}
\DeclareMathOperator{\Gal}{Gal}

\newcommand{\CC}{\mathbb{C}}
\newcommand{\QQ}{\mathbb{Q}}
\newcommand{\FF}{\mathbb{F}}
\newcommand{\ZZ}{\mathbb{Z}}
\newcommand{\GG}{\mathbb{G}}
\newcommand{\TT}{\mathbb{T}}
\newcommand{\HH}{\mathbb{H}}
\newcommand{\NN}{\mathbb{N}}

\newcommand{\cO}{\mathcal{O}}

\newcommand{\Qbar}{{\overline{\QQ}}}
\newcommand{\Fbar}{{\overline{\FF}}}
\newcommand{\Fpbar}{{\overline{\FF_p}}}
\newcommand{\Qpbar}{{\overline{\QQ}_p}}

\newcommand{\res}{\mathrm{res}}
\newcommand{\cores}{\mathrm{cores}}
\newcommand{\torsion}{\mathrm{torsion}}
\newcommand{\modulo}{\,\mathrm{mod}\,}
\newcommand{\mult}{{\rm mult}}
\newcommand{\Coind}{{\rm Coind}}

\newcommand{\End}{{\rm End}}
\newcommand{\Hom}{{\rm Hom}}
\newcommand{\GL}{\mathrm{GL}}
\newcommand{\SL}{\mathrm{SL}}
\newcommand{\Pic}{\mathrm{Pic}}
\newcommand{\cusps}{\mathrm{cusps}}

\newcommand{\cW}{\mathcal{W}}
\newcommand{\fP}{\mathfrak{P}}
\newcommand{\PSL}{\mathrm{PSL}}
\newcommand{\pf}{{\bf Proof. }}
\newcommand{\qed}{\hspace* {.5cm} \hfill $\Box$}
\newcommand{\mat}[4]{
 \left(  \begin{smallmatrix} #1 & #2 \\ #3 & #4 \end{smallmatrix} \right)}
\newcommand{\vect}[2]{
 \left(  \begin{smallmatrix} #1 \\ #2 \end{smallmatrix} \right)}
\newcommand{\dia}[1]{\langle #1 \rangle}
\newcommand{\diap}[1]{\langle #1 \rangle_p}
\newcommand{\diaN}[1]{\langle #1 \rangle_N}
\newcommand{\pair}[2]{\langle #1, #2 \rangle}
\newcommand{\Hpar}{H_{\mathrm{par}}}

\newcommand{\conj}{\mathrm{conj}}
\newcommand{\proj}{\overset{\mathrm{proj}}\twoheadrightarrow}
\newcommand{\Sh}{\mathrm{Sh}}
\newcommand{\Cot}{\mathrm{Cot}}
\newcommand{\Lbar}{{\overline{L}}}
\newcommand{\matz}{\mathrm{Mat}_2(\ZZ)_{\neq 0}}

\begin{document}

\title{On the faithfulness of parabolic cohomology\\ 
as a Hecke module over a finite field}
\author{Gabor Wiese}
\maketitle
 
\begin{abstract}
This article exhibits conditions under which
a certain parabolic group cohomology space
over a finite field $\FF$ is a faithful module for
the Hecke algebra of cuspidal Katz modular forms over
an algebraic closure of $\FF$.
These results can e.g.\ be applied to compute cuspidal Katz
modular forms of weight one with
methods of linear algebra over $\FF$.

MSC Classif.: 11F25 (primary), 11F33, 11F67, 11Y40 (secondary).
\end{abstract}

\section{Introduction}

The space $S_k(\Gamma_1(N),\CC)$ of holomorphic 
modular cusp forms of weight $k \ge 2$ for the group
$\Gamma_1(N)$ for some integer $N \ge 1$ is related to
the parabolic group cohomology of $\Gamma_1(N)$ 
by the so-called Eichler-Shimura isomorphism (see e.g.\ \cite{DiamondIm}, 
Theorem~12.2.2)
\begin{equation}\label{es}
S_k(\Gamma_1(N),\CC) \oplus \overline{S_k(\Gamma_1(N),\CC)} 
  \cong \Hpar^1(\Gamma_1(N),\CC[X,Y]_{k-2}),
\end{equation}
where $\CC[X,Y]_{k-2}$ denotes the $\CC$-vector space of homogeneous
polynomials of degree $k-2$ in two variables with the
standard $\SL_2(\ZZ)$-action.
For the definition of the parabolic group cohomology see 
Section~\ref{SecPar}.
On the modular forms, as well as on the group cohomology
one disposes of natural Hecke operators, which are compatible
with the Eichler-Shimura isomorphism.
An immediate observation is
that $\Hpar^1(\Gamma_1(N),\CC[X,Y]_{k-2})$
is a free module of rank~$2$ for
the complex Hecke algebra $\TT_\CC$ generated by the Hecke operators
in the endomorphism ring of $S_k(\Gamma_1(N),\CC)$. 
This implies, in particular, that
$\Hpar^1(\Gamma_1(N),\CC[X,Y]_{k-2})$ is a faithful $\TT_\CC$-module.

In general there is no analogue of the Eichler-Shimura isomorphism
over finite fields. However, in the present article we show that the
above observation generalises to cuspidal Katz modular forms over finite fields
in certain cases.
\smallskip

Let now $\TT_{\FF_p}$ denote the $\FF_p$-Hecke algebra generated
by the Hecke operators inside the endomorphism ring of the
cuspidal Katz modular forms $S_k(\Gamma_1(N),\FF_p)$ of weight $k\ge 2$
for $\Gamma_1(N)$ over the field~$\FF_p$. 
In view of the observations from the Eichler-Shimura
isomorphism two questions arise naturally.
\begin{enumerate}[(a)]
\itemsep=0cm plus 0pt minus 0pt
\item Is $\Hpar^1(\Gamma_1(N),\FF_p[X,Y]_{k-2})$ a faithful $\TT_{\FF_p}$-module?
\item Is $\Hpar^1(\Gamma_1(N),\FF_p[X,Y]_{k-2})$ free of rank~$2$ as a $\TT_{\FF_p}$-module?
\end{enumerate}
A positive answer to Question~(b) clearly implies a positive answer to~(a).
\smallskip

Using $p$-adic Hodge theory (in particular the article \cite{FaltingsJordan}), 
Edixhoven shows in Theorem~5.2 of \cite{EdixJussieu} that 
Question~(a) has a positive answer in the
weight range $2 \le k \le p-1$.
Emerton, Pollack and Weston were able to deduce from the fundamental work by
Wiles on Fermat's last theorem
that Question~(b) is true locally at primes of the Hecke algebra which
are $p$-ordinary and $p$-distinguished (see \cite{EPW}, Proposition~3.3.1)
for $\ZZ_p$-coefficients and the full
cohomology (i.e.\ not the parabolic subspace), from which the result follows
for $\FF_p$-coefficients and the full cohomology.
We should mention that a positive answer
to Question~(a) with $\ZZ_p$-coefficients does not imply a positive answer over~$\FF_p$.
\smallskip

The main result of the present article is to answer Question~(a) positively
locally at $p$-ordinary primes of $\TT_{\FF_p}$ 
in the weight range $2 \le k \le p+1$ (see Corollary~\ref{corap}).
A similar result for cuspidal modular forms for $\Gamma_1(N)$ with a 
Dirichlet character follows under certain rather weak assumptions
(see Theorem~\ref{heckedelta}).
These results are obtained using techniques that were inspired
by the special case $p=2$ of \cite{EdixJussieu}, Theorem~5.2.
\smallskip

Let us point out the following consequence, which presented
the initial motivation for this study.
If the answer to Question~(a) is yes, we may compute $\TT_{\FF_p}$ via
the Hecke operators on the finite dimensional $\FF_p$-vector space
$\Hpar^1(\Gamma_1(N),\FF_p[X,Y]_{k-2})$.
This can be done explicitly on a computer and only requires 
linear algebra methods.
The knowledge of the Hecke algebra is very interesting. 
One can, for example,
try to verify whether it is a Gorenstein ring. This property is often implied
when the Hecke algebra is isomorphic to a deformation ring (as e.g.\ in \cite{Wiles}).

Techniques from \cite{EdixJussieu}, Section~4, show how cuspidal Katz modular forms
of weight one over $\Fpbar$ can be related to weight~$p$. 
We observe in Proposition~\ref{eigenspace}
that the corresponding local factors of the Hecke algebra 
in weight~$p$ are $p$-ordinary.
Hence, a consequence of our result is that methods from $\FF_p$-linear algebra
can be used for the computation of weight one Katz cusp forms (see Theorem~\ref{wtoneOK}).
The author carried out some such calculations and reported on them
in \cite{app}.
\smallskip

It should be mentioned that if one is not interested in the Hecke
algebra structure, but only in the systems of eigenvalues
of Hecke eigenforms in $S_k(\Gamma_1(N),\Fpbar)$, one can get them
directly from $H^1(\Gamma_1(N),\FF_p[X,Y]_{k-2})$ for all $k \ge 2$
(see e.g.\ \cite{Thesis}, Proposition~4.3.1). 
However, there may be more systems of eigenvalues in the group
cohomology than in $S_k(\Gamma_1(N),\Fpbar)$, but the extra ones can easily
be identified.
\smallskip

This article is based on the use of group cohomology. 
The author studied in his thesis \cite{Thesis} and the article \cite{ArtEins}
relations to modular symbols and a certain cohomology group on the corresponding
modular curve. 
All these agree for instance for any congruence subgroup of $\SL_2(\ZZ)$
if the characteristic of the base field is $p \ge 5$ or in any characteristic for the 
group $\Gamma_1(N)$ with $N \ge 5$.
\smallskip

We add some remarks on Katz modular forms.
As a reference one can use \cite{EdixBoston}, \cite{Katz}
or \cite{DiamondIm}.
For $R=\CC$, the Katz cuspidal modular forms are precisely the holomorphic 
cusp forms.
Under the assumptions $N \ge 5$ and $k \ge 2$
one has $S_k(\Gamma_1(N),\ZZ[1/N]) \otimes R \cong S_k(\Gamma_1(N),R)$ 
for any ring~$R$ in which $N$ is invertible.
In particular, this means that any Katz cusp form over $\Fpbar$ can be
lifted to characteristic zero in the same weight and level.
The previous statement does not hold for weight $k=1$ in general.

Katz cusp forms over~$\Fpbar$ for a prime~$p$ play an important 
r\^ole in a modified version of Serre's conjecture on modular forms
(see \cite{EdixWeight}), as their use, on the one hand, 
gets rid of some character lifting problems, and
on the other hand, allows one to minimise the weights of the modular forms.
In particular, odd, irreducible Galois representations 
$\Gal(\Qbar|\QQ) \to \GL_2(\Fpbar)$
which are unramified at~$p$ are conjectured to correspond to 
Katz eigenforms of weight one.
\smallskip

Finally, we should mention that William Stein has implemented $\FF_p$-modular symbols
and the Hecke operators on them in Magma and Python. One could say that this article
is also about what they compute.

\subsection*{Acknowledgements}

This article has grown out of Chapter~III of my thesis.
I would like to thank my thesis advisor Bas Edixhoven for very helpful suggestions,
and Theo van den Bogaart for sharing his insights in algebraic geometry.

\subsection{Notation}\label{notation}

We denote by $\matz$ the semi-group of integral $2 \times 2$-matrices 
with non-zero determinant. On it one has 
Shimura's main involution ${\mat a b c d}^\iota = \mat d{-b}{-c}a$.
For a ring~$R$ the notation $R[X,Y]_n$ stands for the homogeneous polynomials
of degree~$n$ in the variables $X,Y$.
We let
$$ V_n(R) := \Sym^{n}(R^2) \cong R[X,Y]_n$$ 
which we equip with the natural left $\matz$-semi-group action.
If $V$ is an $R$-module, let $V^\vee$ be the dual $R$-module $\Hom_R(V,R)$.

As usual, the upper half plane is denoted by~$\HH$.
Moreover, we shall use the standard congruence subgroups $\Gamma(N)$,
$\Gamma_1(N)$ and $\Gamma_0(N)$ of $\SL_2(\ZZ)$ for an integer~$N \ge 1$,
which are defined to consist of the matrices $\mat a b c d \in \SL_2(\ZZ)$
reducing modulo~$N$ to $\mat 1 0 0 1$, $\mat 1 * 0 1$ resp.\ $\mat * * 0 *$.

\section{Parabolic group cohomology}\label{SecPar}

In this section we recall the definition of parabolic group
cohomology and present some properties to be used in the
sequel. We shall use standard facts on group cohomology
without special reference. The reader can 
for example consult~\cite{Brown}.

\subsection*{Cohomology of $\PSL_2(\ZZ)$}

The group $\PSL_2(\ZZ)$ is freely generated by the matrix classes of
$\sigma := \mat 0 {-1} 1 0$ and $\tau := \mat 1 {-1} 1 0$.
A consequence is the following special case of the 
{\em Mayer-Vietoris sequence} (see e.g.\ \cite{HiltonStammbach}, p.~221).
For any ring~$R$ and any left $R[\PSL_2(\ZZ)]$-module~$M$ 
the sequence
\begin{multline}\label{mayervietoris}
0 \to M^{\PSL_2(\ZZ)} \to M^{\langle \sigma \rangle} \oplus M^{\langle \tau \rangle} 
\to M \\ \to
H^1(\PSL_2(\ZZ),M) \to H^1(\langle \sigma \rangle, M) \oplus H^1(\langle \tau \rangle, M) 
\to 0
\end{multline}
is exact and for all $i \ge 2$ one has isomorphisms
\begin{equation}\label{mveins}
H^i(\PSL_2(\ZZ),M) \cong H^i(\langle \sigma \rangle, M) \oplus H^i(\langle \tau \rangle, M).
\end{equation}

\begin{cor}\label{corhzweinull}
Let $R$ be a ring and $\Gamma \le \PSL_2(\ZZ)$ be a subgroup of finite index such
that all the orders of all stabiliser groups $\Gamma_x$ for $x \in \HH$ are invertible in~$R$.
Then for all $R[\Gamma]$-modules~$V$ one has
$H^1(\Gamma,V) = M/(M^{\langle \sigma \rangle} + M^{\langle \tau \rangle})$ with $M = \Coind_\Gamma^{\PSL_2(\ZZ)}(V)$
and
$H^i(\Gamma,V) = 0$ for all $i \ge 2$.
\end{cor}

\pf
We first recall that any non-trivially stabilised point $x$ of~$\HH$
is conjugate by an element of $\PSL_2(\ZZ)$ to either $i$ or~$\zeta_3 = e^{2\pi i /3}$, 
whence all non-trivial stabiliser groups are of the form
$g \langle \sigma \rangle g^{-1} \cap \Gamma$
or $g \langle \tau \rangle g^{-1} \cap \Gamma$ for some $g \in \PSL_2(\ZZ)$.
From Mackey's formula (see e.g.\ \cite{ArtEins}) one obtains
$$H^i(\langle \sigma \rangle, \Coind_\Gamma^{\PSL_2(\ZZ)} V) \cong
\prod_{g \in \Gamma \backslash \PSL_2(\ZZ) / \langle \sigma \rangle} 
H^i (g \langle \sigma \rangle g^{-1} \cap \Gamma, V)$$
for all~$i$ and a similar result for~$\tau$. 
Due to the invertibility assumption the result follows from Shapiro's lemma
and Equations (\ref{mayervietoris}) and~(\ref{mveins}).
\qed
\smallskip

The assumptions of the proposition are for instance always satisfied
if $R$ is a field of characteristic not $2$ or $3$. They also
hold for $\Gamma_1(N)$ with $N \ge 4$ over any ring.

\subsection*{Definition of parabolic group cohomology}

Let $R$ be a ring, $\Gamma \le \PSL_2(\ZZ)$ a subgroup of finite index
and $T = \tau \sigma = \mat 1 1 0 1$.
One defines the {\em  parabolic group cohomology group for the $R[\Gamma]$-module~$V$}
as the kernel of the restriction map in
\begin{equation}\label{pardef}
0 \to \Hpar^1(\Gamma, V) \to H^1(\Gamma,V) \xrightarrow{\res} 
\prod_{g \in \Gamma \backslash \PSL_2(\ZZ) / \langle T \rangle}
H^1(\Gamma \cap \langle g T g^{-1} \rangle, V).
\end{equation}
The definition of parabolic cohomology is compatible with Shapiro's lemma,
i.e.\ Equation~(\ref{pardef}) is isomorphic to
\begin{equation}\label{parshapiro}
0 \to \Hpar^1(\PSL_2(\ZZ), M) \to H^1(\PSL_2(\ZZ), M) 
\xrightarrow{\res} H^1(\langle T \rangle, M)
\end{equation}
with $M = \Coind_\Gamma^{\PSL_2(\ZZ)} V = \Hom_{R[\Gamma]}(R[\PSL_2(\ZZ)],V)$,
as one sees using e.g.\ Mackey's formula as in the proof
of Corollary~\ref{corhzweinull}.

\begin{prop}\label{leraypar}
Let $R$ be a ring and $\Gamma \le \PSL_2(\ZZ)$ be a subgroup of finite index such
that all the orders of all stabiliser groups $\Gamma_x$ for $x \in \HH$ are invertible in~$R$.
Then for all $R[\Gamma]$-modules~$V$ the sequence
$$ 0 \to \Hpar^1(\Gamma, V) \to H^1(\Gamma,V) \xrightarrow{\res} 
\prod_{g \in \Gamma \backslash \PSL_2(\ZZ) / \langle T \rangle}
H^1(\Gamma \cap \langle g T g^{-1} \rangle, V) \to V_\Gamma \to 0 $$
is exact.
\end{prop}

\pf
Due to the assumptions we may apply Corollary~\ref{corhzweinull}.
The restriction map in Equation~(\ref{parshapiro}) thus becomes
$$ M/(M^{\langle \sigma \rangle} + M^{\langle \tau \rangle}) \xrightarrow{m \mapsto (1-\sigma)m} M/(1-T)M,$$
since $H^1(\langle T \rangle, M) \cong M/(1-T)M$.
The isomorphism $M \cong (R[\PSL_2(\ZZ)] \otimes_R V)_\Gamma$ allows
one to compute that the cokernel of this map is $V_\Gamma$, the $\Gamma$-coinvariants.
\qed

\subsection*{The module $V_n(R)$}

Let $R$ be a ring.
Recall from Notation~\ref{notation} that we put 
$V_n(R) = \Sym^{n}(R^2) \cong R[X,Y]_n$.

\begin{prop}\label{vkdual}
Suppose that $n!$ is invertible in~$R$. Then there is a perfect pairing
$$V_n(R) \times V_n(R) \to R$$
of $R$-modules, which induces
an isomorphism $V_n(R) \to V_n(R)^\vee$ of $R$-modules 
respecting the $\matz$-action which is given on $V_n(R)^\vee$ by 
$(M.\phi)(w) = \phi(M^\iota w)$ for $M \in \matz$, $\phi \in V_n(R)^\vee$
and $w \in V_n(R)$.
\end{prop}

\pf
One defines the perfect pairing on~$V_n(R)$
by first constructing a perfect pairing on $R^2$, which we consider
as column vectors. One sets
$$ R^2 \times R^2 \to R, \;\;  \pair v w := 
\det ( v | w )= v_1 w_2 - v_2 w_1.$$
If $M$ is a matrix in $\matz$, one checks easily that
$\pair  {M v}{w} =  \pair v {M^\iota w}$. 
This pairing extends naturally to a pairing on the $n$-th tensor power of~$R^2$.
Due to the assumption on the invertibility of~$n!$, 
we may view $\Sym^n(R^2)$ as a submodule in the $n$-th tensor power, and hence
obtain the desired pairing and the isomorphism of the statement.
\qed

\begin{lem}\label{lemvknull}
Let $n \ge 1$ be an integer, $t = \mat 1 N 0 1$ and
$t' = \mat 1 0 N 1$.
If $n! N$ is not a zero divisor in~$R$, then
for the $t$-invariants we have
$V_n(R)^{\langle t \rangle} = \langle X^n \rangle$
and for the $t'$-invariants
$V_n(R)^{\langle t' \rangle} = \langle Y^n \rangle$.
If $n! N$ is invertible in~$R$, then the coinvariants are
given by
$V_n(R)_{\langle t \rangle} = V_n(R)/\langle Y^n, XY^{n-1}, \dots, X^{n-1}Y \rangle$
respectively
$V_n(R)_{\langle t' \rangle} = V_n(R)/\langle X^n, X^{n-1}Y, \dots, XY^{n-1} \rangle$.
\end{lem}

\pf
The action of~$t$ is $t.(X^{n-i} Y^i) = X^{n-i} (NX + Y)^i$ and consequently
$(t-1). (X^{n-i} Y^i) = \sum_{j=0}^{i-1} r_{i,j} X^{n-j}Y^j$
with $r_{i,j} =N^{i-j} \vect{i}{j}$, which is not a zero divisor, 
respectively invertible, by assumption.
For $x = \sum_{i=0}^n a_i X^{n-i}Y^i$ we have
$ (t-1).x = \sum_{j=0}^{n-1} X^{n-j}Y^j (\sum_{i=j+1}^n a_i r_{i,j}).$
If $(t-1).x = 0$, we conclude for $j = n-1$ that $a_n = 0$. 
Next, for $j = n-2$ it follows that $a_{n-1} = 0$,
and so on, until $a_1 = 0$. This proves the statement on
the $t$-invariants. The one on the $t'$-invariants follows from symmetry. 
The claims on the coinvariants are proved in a very similar and
straightforward way. 
\qed

\begin{prop}\label{vknull}
Let $n\ge 1$ be an integer. 
\begin{enumerate}[(a)]
\item If $n! N$ is not a zero divisor in~$R$, 
then the $R$-module of $\Gamma(N)$-invariants $V_n(R)^{\Gamma(N)}$ is zero. 
\item If $n! N$ is invertible in~$R$,
then the $R$-module of $\Gamma(N)$-coinvariants $V_n(R)_{\Gamma(N)}$ is zero. 
\item Suppose that $\Gamma$ is a subgroup of $\SL_2(\ZZ)$ such
that reduction modulo~$p$ defines a surjection $\Gamma \twoheadrightarrow \SL_2(\FF_p)$
(e.g.\ $\Gamma(N)$, $\Gamma_1(N)$, $\Gamma_0(N)$ for $p \nmid N$). 
Suppose moreover that $1 \le n \le p$ if $p>2$, and $n=1$ if $p=2$.
Then one has $V_n(\FF_p)^{\Gamma} = 0 = V_n(\FF_p)_\Gamma.$
\end{enumerate}
\end{prop}

\pf
As $\Gamma(N)$ contains the matrices $t$ and $t'$, Lemma~\ref{lemvknull}
already finishes Parts~(a) and~(b). 
The only part of~(c) that is not yet covered is when the degree is $n=p>2$.
Since $V_p(\FF_p)$ is naturally isomorphic to~$U_1$,
Proposition~\ref{propudvd} gives the exact sequence of $\Gamma$-modules
$0 \to V_1(\FF_p) \to V_p(\FF_p) \to V_{p-2}(\FF_p) \to 0.$
It suffices to take invariants respectively coinvariants to obtain the result.
\qed

\subsection*{Torsion-freeness and base change properties}

Herremans has computed a torsion-freeness result like the following proposition
in \cite{Herremans}, Proposition~9. 
Here we give a short and conceptual proof of a slightly more general statement.
The way of approach was suggested by Bas Edixhoven.

\begin{prop}\label{torsion}
Assume that $R$ is an integral domain of characteristic~$0$
such that $R/pR \cong \FF_p$ for a prime~$p$.
Let $N \ge 1$ and $k \ge 2$ be integers and let $\Gamma \le \SL_2(\ZZ)$
be a subgroup containing $\Gamma(N)$ but not $-1$ such that
the orders of the stabiliser subgroups $\Gamma_x$ for $x \in \HH$
have order coprime to~$p$.
Then the following statements hold:
\begin{enumerate}[(a)]
\item $H^1(\Gamma,V_{k-2}(R)) \otimes_R \FF_p \cong H^1(\Gamma,V_{k-2}(\FF_p))$.
\item If $k=2$, then $H^1(\Gamma,V_{k-2}(R))[p]=0$.
If $k\ge 3$, then $H^1(\Gamma,V_{k-2}(R))[p]=V_{k-2}(\FF_p)^\Gamma$.
In particular, if $p \nmid N$, then $H^1(\Gamma,V_{k-2}(R))[p]=0$ 
for all $k \in \{2, \dots, p+2\}$.
\item If $k=2$, or if $k \in \{3, \dots, p+2\}$ and $p \nmid N$, then
$\Hpar^1(\Gamma,V_{k-2}(R)) \otimes_R \FF_p \cong \Hpar^1(\Gamma,V_{k-2}(\FF_p))$.
\end{enumerate}
\end{prop}

\pf
Let us first notice that the sequence
$$ 0 \to V_{k-2}(R) \xrightarrow{\cdot p} V_{k-2}(R) \to 
V_{k-2}(\FF_p) \to 0$$
of $R[\Gamma]$-modules is exact. 
The associated long exact sequence 
gives rise to the short exact sequence
$$ 0 \to H^i(\Gamma,V_{k-2}(R)) \otimes \FF_p 
     \to H^i(\Gamma, V_{k-2}(\FF_p)) 
     \to H^{i+1}(\Gamma,V_{k-2}(R))[p] \to 0$$
for every $i \ge 0$.
Exploiting this sequence for $i=1$ immediately yields Part~(a), since
any $H^2$ of $\Gamma$ is zero by Corollary~\ref{corhzweinull}.
Part~(b) is a direct consequence of the case $i=0$ and
Proposition~\ref{vknull}.

We have the exact commutative diagram
$$ \xymatrix@=.25cm{
0 \ar@{->}[r] & H^1(\Gamma,V_{k-2}(R)) \ar@{->}[d] \ar@{->}[r]^{\cdot p} &  
H^1(\Gamma,V_{k-2}(R)) \ar@{->}[d] \ar@{->}[r]  & 
H^1(\Gamma,V_{k-2}(\FF_p)) \ar@{->}[d] \ar@{->}[r] & 0 \\
0 \ar@{->}[r] & \prod_g H^1(D_g,V_{k-2}(R)) \ar@{->}[d] \ar@{->}[r]^{\cdot p} &  
\prod_g H^1(D_g,V_{k-2}(R)) \ar@{->}[d] \ar@{->}[r]  & 
\prod_g H^1(D_g,V_{k-2}(\FF_p)) \ar@{->}[r] & 0 \\
& (V_{k-2}(R))_\Gamma \ar@{->}[r]^{\cdot p} \ar@{->}[d] & (V_{k-2}(R))_\Gamma \ar@{->}[d]\\
& 0 & 0 }$$%
where the products are taken over 
$g \in \Gamma \backslash \PSL_2(\ZZ) / \langle T \rangle$,
and $D_g = \Gamma \cap \langle gTg^{-1}\rangle$.
The exactness of the first row is the contents of Parts~(a) and~(b).
That the columns are exact follows from Proposition~\ref{leraypar}.
The zero on the right of the second row is due to the fact that
$D_g$ is free on one generator.
That generator is of the form $g \mat 1 r 0 1 g^{-1}$ with $r \mid N$, 
so that $r$ is invertible in~$\FF_p$. The zero on the left is 
trivial for $k=2$ and for $3 \le k \le p+2$ it is
a consequence of Lemma~\ref{lemvknull}.
Part~(c) now follows from the snake lemma
and Proposition~\ref{vknull}, which implies that the bottom
map is an injection.
\qed

\section {Hecke action}\label{SecHecke}

Hecke operators conceptually come from Hecke correspondences
on modular curves respectively modular stacks. They are best described on 
the moduli interpretations (see e.g.\ \cite{DiamondIm}, 3.2 and 7.3).
All Hecke operators that we will encounter in this article
arise like this.
This section presents Hecke operators on group
cohomology and the principal result is the behaviour of the Hecke
operators with respect to Shapiro's lemma. That result
was obtained by Ash and Stevens (\cite{AshStevens}, Lemma~2.2). Here, however,
we avoid their rather heavy language of weakly compatible Hecke
pairs. Instead, the description of Hecke operators 
on group cohomology is used which comes directly from
the Hecke correspondences (formally one has to work on the modular stacks
with locally constant coefficients, in case of non-trivially 
stabilised points). For the description we follow \cite{DiamondIm}, 12.4.

\subsection* {Hecke operators on group cohomology}

Let $R$ be a ring, $\alpha \in \matz$ and $\Gamma \le \PSL_2(\ZZ)$ 
be a subgroup containing some $\Gamma(N)$.
We use the notations $\Gamma_\alpha := \Gamma \cap \alpha^{-1} \Gamma \alpha$
and $\Gamma^\alpha := \Gamma \cap \alpha \Gamma \alpha^{-1}$,
where we consider $\alpha^{-1}$ as an element of $\GL_2(\QQ)$.
Both groups are commensurable with~$\Gamma$.

Suppose that $V$ is an $R$-module with a $\matz$-semi-group action
which restricts to an action by~$\Gamma$.
The {\em Hecke operator} $T_\alpha$ acting on group
cohomology is the composite
$$ H^1(\Gamma,V) \xrightarrow{\res} H^1(\Gamma^\alpha, V)
 \xrightarrow{\conj_\alpha} H^1(\Gamma_\alpha, V)
 \xrightarrow{\cores} H^1(\Gamma,V).$$
The first map is the usual {\em restriction}, and the third one is the
so-called {\em corestriction}, which one also finds in the
literature under the name {\em transfer}. 
We explicitly describe the second map on non-homogeneous cocycles
(cf.\ \cite{DiamondIm}, p.~116):
$$ \conj_\alpha: H^1(\Gamma^\alpha, V) \to H^1(\Gamma_\alpha, V), \;\;
c \mapsto \big( g_\alpha \mapsto \alpha^{\iota}.c(\alpha g_\alpha \alpha^{-1}) \big).$$
There is a similar description on the parabolic subspace 
and the two are compatible.
The following formula can also be found in \cite{DiamondIm}, p.~116, 
and \cite{Shimura}, Section~8.3.

\begin{prop}\label{shdesc}
Suppose that $\Gamma \alpha \Gamma = \bigcup_{i=1}^n \Gamma \delta_i$ is a disjoint
union. Then the Hecke operator $T_\alpha$ acts on $H^1(\Gamma,V)$ 
and $\Hpar^1(\Gamma,V)$ by
sending the non-homogeneous cocyle $c$ to $T_\alpha c$ defined by
$$ (T_\alpha c)(g) = \sum_{i=1}^n \delta_i^\iota c(\delta_i g \delta_{j(i)}^{-1})$$
for $g \in \Gamma$. Here $j(i)$ is the index such that $\delta_i g \delta_{j(i)}^{-1} \in \Gamma$.
\end{prop}

\pf
We only have to describe the corestriction explicitly. For that we notice that one has
$\Gamma = \bigcup_{i=1}^n \Gamma_\alpha g_i$
with $\alpha g_i = \delta_i$. Furthermore the corestriction of a non-homogeneous
cocycle $u \in H^1(\Gamma_\alpha, V)$ is the cocycle $\cores(u)$ uniquely
given by 
$$\cores(u)(g) = \sum_{i=1}^n g_i^{-1} u(g_i g g_{j(i)}^{-1})$$
for $g \in \Gamma$. Combining with the explicit description of the map $\conj_\alpha$
yields the result.
\qed
\smallskip

Suppose now that $\Gamma = \Gamma_1(N)$ (resp.\ $\Gamma = \Gamma_0(N)$).
For a positive integer~$n$, the {\em Hecke operator} $T_n$ is $\sum_\alpha T_\alpha$,
where the sum runs through a system of representatives of the double
cosets $\Gamma \backslash \Delta^n / \Gamma$ for
the set $\Delta^n$ of matrices $\mat abcd \in \matz$ of
determinant~$n$ such that $\mat abcd \equiv \mat 1*0* \modulo N$ 
(resp.\ $\mat abcd \equiv \mat **0* \modulo N$).
For a prime $p$ one has $T_p = T_\alpha$ with $\alpha = \mat 100p$.
If $\Gamma = \Gamma_1(N)$
and the integer $d$ is coprime to~$N$,
the {\em diamond operator}
$\dia{d}$ is $T_\alpha$ for any matrix $\alpha \in \SL_2(\ZZ)$ whose reduction
modulo $N$ is $\mat {d^{-1}} 0 0 d$.
The diamond operator gives a group action by~$(\ZZ/N\ZZ)^*$
(with $-1$ acting trivially).
If the level is $NM$ with $(N,M)=1$, then we can separate the 
diamond operator into two parts
$\dia d = \langle d \rangle_M \times \diaN d,$
corresponding to $\ZZ/NM\ZZ \cong \ZZ/M\ZZ \times \ZZ/N\ZZ$.
The Hecke and diamond operators satisfy the ``usual'' Euler product
and one has the formulae
$T_n T_m = T_{nm}$ for any pair
of coprime integers $n,m$ and 
$T_{p^{r+1}} = T_{p^r}T_p - p^{k-1} \langle p \rangle T_{p^{r-1}}$
if $p \nmid N$, and $T_{p^{r+1}} = T_{p^r}T_p$ if $p \mid N$
for the action on $H^1(\Gamma_1(N),V_{k-2})$.

\subsection*{Hecke operators and Shapiro's lemma}

Let $N,M$ be coprime positive integers.
In order to test compatibility of the Hecke operators
with Shapiro's lemma we need to extend the $\Gamma_1(N)$-action on 
$\Coind_{\Gamma_1(NM)}^{\Gamma_1(N)}(V)$ to a $\matz$-semi-group action
``in the right way''. For that we define the $R$-module $\cW(M,V)$ as
$$\{f \in \Hom_R (R[(\ZZ/M\ZZ)^2],V) \;|\; f((u,v)) = 0 \;\forall (u,v) \text{ s.t. }
\langle u,v \rangle \neq \ZZ/M\ZZ\}.$$
It carries the left $\matz$-semi-group action
$(g.f)((u,v)) = g f((u,v)g)$ for $f \in \cW(M,V)$, $g \in \matz$ and
$(u,v) \in (\ZZ/M\ZZ)^2$.

\begin{lem}\label{lemsh}
The homomorphism
$$\cW(M,V) \to \Hom_{R[\Gamma_1(NM)]}(R[\Gamma_1(N)],V), \;\;\;
   f \mapsto \big(g \mapsto (g.f)((0,1))\big)$$
is an isomorphism of left $\Gamma_1(N)$-modules (for the restricted action on $\cW(M,V)$).
In particular, $\cW(M,V)$ is isomorphic to $\Coind_{\Gamma_1(NM)}^{\Gamma_1(N)}(V)$
as a left $\Gamma_1(N)$-module.
\end{lem}

\pf
As $N$ and $M$ are coprime, reduction modulo~$M$ defines a surjection from
$\Gamma_1(N)$ onto $\SL_2(\ZZ/M\ZZ)$. This implies that the map
$$ \Gamma_1(NM) \backslash \Gamma_1(N) \xrightarrow {A \mapsto (0,1)A \text{ mod } M}
(\ZZ/M\ZZ)^2$$
is injective, and its image is the set of the $(u,v)$ with $\ZZ/M\ZZ = \langle u,v\rangle$.
As the $\Gamma_1(NM)$-action on~$V$ is the restriction of a $\Gamma_1(N)$-action,
the coinduced module can be identified with 
$\Hom_R(R[\Gamma_1(NM)\backslash \Gamma_1(N)],V)$.
From this the claimed isomorphism follows directly.
\qed

\begin{lem}\label{cosetdec}
Let $N,n$ be positive integers.
We set 
$$\Delta_1(N)^n = \{ \mat abcd \in \matz | \det \mat abcd = n, \mat abcd \equiv \mat 1*0* 
\modulo N\}.$$
There is the decomposition
$$ \Delta_1(N)^n = \bigcup_a \bigcup_b \Gamma_1(N) \sigma_a \mat a b 0 d$$
where $a$ runs through the integers such that $a>0$, $(a,N)=1$, $ad=n$ and $b$ through a system
of representatives of $\ZZ/d\ZZ$. Here $\sigma_a \in \SL_2(\ZZ)$ is a matrix reducing
to~$\mat {a^{-1}} 0 0 a$ modulo~$N$.
\end{lem}

\pf
This is \cite{Shimura}, Proposition 3.36.
\qed
\smallskip

The Shapiro map is the isomorphism on cohomology groups
$$ \Sh: H^1(\Gamma_1(N),\cW(M,V)) \to H^1(\Gamma_1(NM),V)$$
which is induced by the homomorphism $\cW(M,V) \to V$,
sending $f$ to $f((0,1))$.

\begin{prop}\label{shapiro}
For all integers $n,d \ge 1$ with $(d,N)=1$ we have
$$ T_n \circ \Sh = \Sh \circ T_n \;\;\; \text{ and } \;\;\;
\langle d \rangle_N \circ \Sh = \Sh \circ \langle d \rangle_N.$$
\end{prop}

\pf
We first prove the statement for~$T_n$.
For every integer $a>0$ dividing~$n$ such that $(a,N)=1$ 
we choose a matrix~$\sigma_a$ such that it reduces to 
$\mat {a^{-1}} 0 0 a$ modulo~$N$.
If $(a,M)=1$, then we also impose that $\sigma_a$
reduces to $\mat {a^{-1}} 0 0 a$ modulo~$M$. 
If $(a,M)\neq 1$, then we want $\sigma_a \equiv \mat 1 0 0 1$ modulo~$M$.
A simple calculation shows that $(0,1) \big( \sigma_a \mat a b 0 d \big)^\iota$
is congruent to $(0,1)$ modulo~$M$ if $(a,M)=1$ respectively to $(0,a)$ if $(a,M)\neq 1$. 

A set of explicit coset representatives of 
$\Gamma_1(NM) \backslash \Delta_1(NM)^n$
is given by Lemma~\ref{cosetdec} as a subset 
of coset representatives of 
$\Gamma_1(N) \backslash \Delta_1(N)^n$,
namely of those with $(a,M)=1$.

Let now $c \in H^1(\Gamma_1(N),\cW(M,V))$ be a cocycle. 
Then by Proposition~\ref{shdesc} and the definition of the $\matz$-action on
$\cW(M,V)$ we have for $g \in \Gamma_1(NM)$
$$ (\Sh (T_n c))(g) = 
\sum_\delta \delta^\iota \big(c(\delta g \widetilde{\delta}^{-1})((0,1)\delta^\iota)\big),$$
where the sum runs over the above coset representatives for 
$\Gamma_1(N)\backslash \Delta_1(N)^n$ and $\widetilde{\delta}$
is chosen among these representatives such that 
$\delta g \widetilde{\delta}^{-1} \in \Gamma_1(NM)$.
Moreover, we have
$$ (T_n(\Sh (c)))(g) = \sum_\delta \delta^\iota(c(\delta g \widetilde{\delta}^{-1})((0,1))),$$
where now the sum only runs through the subset described above.
By what we have remarked right above $(0,1)\delta^\iota$ is congruent 
to $(0,1)$ modulo $M$ if and only if
$(a,M)=1$. If $(a,M)\neq 1$, then  $\langle a \rangle \neq \ZZ/M\ZZ$,
but we have $(0,1)\delta^\iota \equiv (0,a) \modulo N$.
This proves the compatibility for~$T_n$.

Since $\diaN d$ only depends on $d$ modulo~$N$, we may suppose
that $d$ is congruent to $1$ modulo~$M$. We choose $\alpha = \sigma_d$
such that it reduces to $\mat {d^{-1}}00d$ modulo~$N$ 
and to $\mat 1001$ modulo~$M$.
For $c \in H^1(\Gamma_1(N),\cW(M,V))$ and $g \in \Gamma_1(NM)$
the formula
$$ (\Sh (\diaN d c))(g) = 
\alpha^\iota \big(c(\alpha g \alpha^{-1})((0,1)\alpha^\iota)\big)
= (\diaN d (\Sh (c)))(g)$$
follows.
\qed

\begin{prop}\label{dia}
For $(n,M)=1$ we define the $R[\matz]$-isomorphism
$$ \mult_n: \cW(M,V) \to \cW(M,V),\;\; f \mapsto ((u,v) \mapsto f\big((nu,nv)\big)).$$
Then we have
$$ \langle n \rangle_M \circ \Sh = \Sh \circ \mult_n.$$
\end{prop}

\pf
Let $\sigma \in \SL_2(\ZZ)$ be a matrix
reducing to $\mat {n^{-1}} 0 0 n$ modulo~$M$ 
and to $\mat 1 0 0 1$ modulo~$N$. This means in particular that 
$\sigma \in \Gamma_1(N)$. Hence, for a cocycle $c \in H^1(\Gamma_1(N), \cW(M,V))$
we have 
$$\sigma^{-1} c(\sigma g \sigma^{-1}) = c(g) + (g-1)c(\sigma^{-1}),$$ 
so that the equality $c(\sigma g \sigma^{-1}) = \sigma c(g)$
holds in $H^1(\Gamma_1(N), \cW(M,V))$.

We can now check the claim.
First we have
$$(\langle n \rangle_M \circ \Sh)(c)(g) =
 \sigma^\iota \big( (\sigma.c(g)) ((0,1)) \big) = c(g)((0,1)\sigma).$$
This agrees with $(\Sh \circ \mult_n)(c)(g) = c(g)((0,n))$.
\qed

\section {Reduction to weight $2$ for parabolic group cohomology}

This section deals with the reduction to weight $2$ for parabolic group
cohomology, by which we mean that the parabolic group
cohomology over $\FF_p$ for $\Gamma_1(N)$ with $p \nmid N$
and weight $3 \le k \le p+1$ is related to the weight two parabolic
group cohomology of $\Gamma_1(Np)$. 
Hence, the parabolic group cohomology shows a
behaviour similar to that of modular forms (cf.\ Proposition~\ref{propgross}).
Roughly speaking, the reduction to weight~$2$ on group cohomology comes from the
decomposition of $\Coind_{\Gamma_1(Np)}^{\Gamma_1(N)}(\FF_p)$ into
simple $\FF_p[\SL_2(\FF_p)]$-modules. Those are precisely the
$V_d(\FF_p)$ for $0 \le d \le p-1$.

The contents of this section is already partly present in~\cite{AshStevens}.
However, in that paper the parabolic subspace is not treated.

\subsection* {Decomposition of $\cW(p,\FF_p)$ as $\FF_p[\matz]$-module}

We now relate the $\FF_p[\matz]$-modules
$\cW(p,\FF_p)$ and $V_d(\FF_p)$ for $0 \le d \le p-1$.
It is easy to check that evaluation of polynomials on $\FF_p^2$ 
induces an isomorphism of $\FF_p[\matz]$-modules
\begin{equation}
 \FF_p[X,Y]/(X^p-X, Y^p-Y) \cong \FF_p^{\FF_p^2}.
\end{equation}
We can thus identify $\cW(p,\FF_p)$ with
$\{f \in \FF_p[X,Y]/(X^p-X, Y^p-Y) \;|\; f((0,0)) = 0\}$.
Let $U_d(\FF_p)$ be the subspace consisting of polynomial classes 
of degree~$d$, i.e.\ those that satisfy
$f(lx,ly) = l^d f(x,y)$ for all $l \in \FF_p^*$. Note that the degree is naturally
defined modulo~$p-1$. It is clear that the natural $\matz$-action respects the degree.
We remark that the dimension of any $U_d$ is $p+1$.
By collecting monomials we obtain
\begin{equation}
\cW(p,\FF_p) = \bigoplus_{d=0}^{p-2} U_d(\FF_p).
\end{equation}
Furthermore, we dispose of the perfect bilinear pairing
$$ \cW(p,\FF_p) \times \cW(p,\FF_p) \to \FF_p, \;\;  
   \pair fg = \sum_{(a,b) \in \FF_p^2} f(a,b)g(a,b).$$
It is elementary to check that $U_d(\FF_p)$ pairs to zero 
with $U_e(\FF_p)$ if $(p-1) \nmid (d+e)$. 
Hence, the restricted pairing $U_d(\FF_p) \times  U_{p-1-d}(\FF_p) \to \FF_p$
is perfect for $0 \le d \le p-1$, 
as the dimensions of $U_{p-1-d}(\FF_p)$ and $U_d(\FF_p)$ are equal.
Furthermore, $\FF_p[X,Y]_d$ pairs to zero with $\FF_p[X,Y]_{p-1-d}$.
Consequently the induced pairing 
$U_d(\FF_p)/V_d(\FF_p) \times V_{p-1-d}(\FF_p) \to \FF_p$
is perfect.
Composing with the map from Proposition~\ref{vkdual}, we obtain an isomorphism
$U_d(\FF_p) / V_d(\FF_p) \to V_{p-1-d}(\FF_p)$.

\begin{prop}\label{propudvd}
Let $p$ be a prime and $d$ and integer with $0 \le d \le p-1$.
The preceding construction gives the exact sequence
$$ 0 \to V_d(\FF_p) \to U_d(\FF_p) \to V_{p-1-d}(\FF_p) \to 0.$$
Moreover, the diagram
$$ \xymatrix@=.5cm{
0 \ar@{->}[r] & V_d(\FF_p) \ar@{->}[d]_{M.} \ar@{->}[r]& U_d(\FF_p)  \ar@{->}[d]^{M.} \ar@{->}[r]
& V_{p-1-d}(\FF_p) \ar@{->}[d]^{\det(M)^d M.} \ar@{->}[r] & 0\\
0 \ar@{->}[r] & V_d(\FF_p) \ar@{->}[r]& U_d(\FF_p)  \ar@{->}[r]
& V_{p-1-d}(\FF_p)  \ar@{->}[r] & 0 }$$
commutes for all $M \in \matz$ whose reduction modulo~$p$ is invertible.
If $d > 0$, then the diagram is also commutative 
for $M = {\mat 1 0 0 p}^\iota = \mat p 0 0 1$, which means in particular
that the right hand side vertical arrow is zero.
\end{prop}

\pf
Only the commutativity need be checked. If $M$ is invertible modulo~$p$,
one immediately sees that the pairing on $\cW(p,\FF_p)$ is invariant under~$M$,
i.e.\ $\langle Mf, Mg \rangle = \langle f, g \rangle$ for
all $f,g \in \cW(p,\FF_p)$.
The commutativity is then clear from Proposition~\ref{vkdual}.
In order to treat $M = \mat p 0 0 1$ one considers the
basis of $U_d(\FF_p)$ given by the monomials of degree~$d$, which correspond to
the embedding of~$V_d(\FF_p)$, together with the monomials 
$X^iY^{p-1+d-i}$ for $d \le i \le p-1$.
As the latter monomials all contain at least one factor of~$X$ by assumption, 
they are killed by applying the matrix.
\qed

\subsection*{Reduction to weight $2$}

We introduce the following notation. Let $M$ be any $\FF_p$-vector space
on which the Hecke operators $T_l$ and the $p$-part of the diamond operators $\diap \cdot$ act.
Suppose $d > 0$.
By $M[d]$ we mean $M$ with the action of the Hecke operator $T_l$ ``twisted''  
to be $l^{d} T_l$ (in particular $T_p$ acts as zero). 
Furthermore, by $M(d)$ we denote the subspace
on which $\diap l$ acts as $l^d$.

The non-parabolic part of the following proposition is also
\cite{AshStevens}, Theorem~3.4.

\begin{prop}\label{heckepar}
Let $p$ be a prime, $N \ge 4$ and $0 < d \le p-1$ integers such that $p \nmid N$.
We have isomorphisms respecting the Hecke operators
$$H^1 (\Gamma_1(Np), \FF_p)(d) \cong H^1(\Gamma_1(N),U_d(\FF_p)) \;\;\;\text{ and}$$
$$\Hpar^1 (\Gamma_1(Np), \FF_p)(d) \cong \Hpar^1(\Gamma_1(N),U_d(\FF_p)).$$
Moreover, the sequences
\begin{multline*}
0 \to H^1(\Gamma_1(N),V_d(\FF_p)) \to H^1(\Gamma_1(N),U_d(\FF_p))\\ \to
H^1(\Gamma_1(N),V_{p-1-d}(\FF_p))[d] \to 0
\end{multline*}
and
\begin{multline*}
0 \to \Hpar^1(\Gamma_1(N),V_d(\FF_p)) \to \Hpar^1(\Gamma_1(N),U_d(\FF_p))\\
\to \Hpar^1(\Gamma_1(N),V_{p-1-d}(\FF_p))[d] \to 0
\end{multline*}
are exact and respect the Hecke operators.
\end{prop}

\pf
The first statement follows from Propositions~\ref{shapiro} and~\ref{dia}
together with the definition of $U_d(\FF_p)$.
The first exact sequence is part of the long exact sequence
associated to the short exact sequence of Proposition~\ref{propudvd}.
The commutative diagram in that proposition
gives the twisting of the Hecke action in the exact sequences.
This follows directly from the description of the Hecke operator
on group cohomology.

For $d=p-1$ we have $U_{p-1}(\FF_p) = V_0(\FF_p) \oplus V_{p-1}(\FF_p)$,
from which the statements follow. So we now assume $d < p-1$,
in particular $p \neq 2$.
For the top sequence we only need to check that it is exact on the left
and on the right. By Proposition~\ref{vknull} we have
$H^0(\Gamma_1(N),V_{p-1-d}(\FF_p)) = 0$.
The $H^2$-terms are trivial by Corollary~\ref{corhzweinull}.

The exactness of the second sequence follows from the snake lemma,
once we have established the exactness of
\begin{multline*}
 0\to \prod_{c \text{ cusps}}H^1(D_c,V_d(\FF_p)) \to 
        \prod_{c \text{ cusps}}H^1(D_c,U_d(\FF_p)) \\ \to 
        \prod_{c \text{ cusps}}H^1(D_c,V_{p-1-d}(\FF_p)) \to 0,
\end{multline*}
where $D_c$ is the stabiliser group of the cusp $c = g\infty$
with $g \in \PSL_2(\ZZ)$.
Hence, $D_c = g \langle T \rangle g^{-1} \cap \Gamma_1(N)$. This
group is infinite cyclic generated by 
$g \mat 1 r 0 1 g^{-1}$ for some $r \in \ZZ$ dividing~$N$ (see also the proof
of Proposition~\ref{torsion}).
Hence, we have $H^2(D_c,V_d(\FF_p))=0$.
We claim that the sequence
\begin{multline*}
 0\to \prod_{c \text{ cusps}}H^0(D_c,V_d(\FF_p)) \to 
        \prod_{c \text{ cusps}}H^0(D_c,U_d(\FF_p)) \\ \to 
        \prod_{c \text{ cusps}}H^0(D_c,V_{p-1-d}(\FF_p)) \to 0
\end{multline*}
is exact. Lemma~\ref{lemvknull} implies
that $H^0(D_c,V_d(\FF_p))$ and $H^0(D_c,V_{p-1-d}(\FF_p))$ are
one-dimensional.
In order to finish the proof, it thus suffices to prove that $H^0(D_c,U_d(\FF_p))$
is (at least) of dimension~$2$. The elements $X^d \in U_d(\FF_p)$ 
and $Y^d(1-X^{p-1}) \in U_d(\FF_p)$ are invariant under~$T$. Indeed, 
\begin{multline*}
T.Y^d(1-X^{p-1}) = (X+Y)^d(1-X^{p-1}) \\
 = Y^d(1-X^{p-1}) + \sum_{i=1}^d \vect d i Y^{d-i} X^i(1-X^{p-1}) = Y^d(1-X^{p-1}),
\end{multline*}
as in $U_d(\FF_p)$ we have $X^i(1-X^{p-1}) = X^{i-1} (X-X^p) = 0$ for $i > 0$.
\qed

\section{Modular forms of weight~$2$ and level $Np$}

This section recalls the reduction to weight $2$
for modular forms and establishes
a link between parabolic group cohomology and modular forms via
Jacobians of modular curves.

\subsection*{Reduction to weight $2$  for modular forms}

We recall some work of Serre as explained in~\cite{Gross}, Sections 7 and~8,
cf.\ also \cite{EdixWeight}, Section~6.

Let us now introduce notation that is used 
throughout the sequel of this article.
We consider the modular curve $X_1(Np)$ over $\QQ_p(\zeta_p)$
for a prime~$p>2$ not dividing~$N \ge 5$.
It has a regular stable model~$X$ over the ring $\ZZ_p[\zeta_p]$, see
e.g.\ \cite{KM}.
Let $J$ denote the N\'eron model over $\ZZ_p[\zeta_p]$ of $J_1(Np)$,
the Jacobian of $X_1(Np)$ over $\QQ_p(\zeta_p)$.
We let, following \cite{Gross}, Section~8,
$$ L = H^0(X, \Omega_{X/\ZZ_p[\zeta_p]}),$$
where $\Omega_{X/\ZZ_p[\zeta_p]}$ is the {\em dualising sheaf of $X$}
of \cite{DR}, Section~I.2. By \cite{Gross}, Equation~8.2, we have
for the special fibre $X_{\FF_p}$ that
$$ \Lbar := H^0(X_{\FF_p}, \Omega_{X_{\FF_p}/\FF_p}) = L \otimes_{\ZZ_p[\zeta_p]} \FF_p.$$

On $L$ and $\Lbar$ the $p$-part $\diap \cdot$ 
of the diamond operator acts.
The principal result on $\Lbar$ that we will need is the following
(see Propositions 8.13 and~8.18 of \cite{Gross}). The
notation is as in Proposition~\ref{heckepar}.

\begin{prop}[Serre]\label{propgross}
Assume $3 \le k \le p$, $N\ge 5$ and $p\nmid N$.
Then there is an isomorphism of $\FF_p$-vector spaces
$$ \Lbar(k-2) \cong S_k(\Gamma_1(N),\FF_p) \oplus S_{p+3-k}(\Gamma_1(N),\FF_p)[k-2]$$
respecting the Hecke operators. Moreover, the sequence of Hecke modules
$$ 0 \to S_2(\Gamma_1(N),\FF_p)[p-1] \to \Lbar(p-1) \to S_{p+1}(\Gamma_1(N),\FF_p) \to 0$$
is exact.
\end{prop}

\subsection* {Parabolic cohomology and the $p$-torsion of the Jacobian}

In order to compare Hecke algebras of cusp forms with
those of parabolic group cohomology in characteristic~$p$, we generalise the
strategy of the second part of the proof of \cite{EdixJussieu}, Theorem~5.2.
Hence, we wish to bring the Jacobian into the play, since it will enable us
to pass from characteristic zero geometry to characteristic~$p$.

\begin{lem}\label{cotlem}
There are isomorphisms
$$\Lbar \cong \Cot_0(J^0_{\FF_p}) \cong \Cot_0(J^0_{\FF_p}[p])$$
respecting the Hecke operators.
\end{lem}

\pf
The first isomorphism is e.g.\ \cite{EdixWeight}, Equation~6.7.2.
The second one follows from the fact that multiplication by~$p$ on
$J_{\FF_p}^0$ induces multiplication by~$p$ on the tangent space at~$0$, 
which is the zero map. Hence, the tangent space at~$0$ of $J_{\FF_p}^0[p]$
is equal to the one of~$J_{\FF_p}^0$.
\qed
\smallskip

To establish an explicit link between parabolic cohomology and modular forms,
we identify the parabolic cohomology group for $\Gamma_1(N)$ with 
$\FF_p$-coefficients as the $p$-torsion of the Jacobian of the corresponding 
modular curve. Here we may view the Jacobian as a complex abelian variety.
As in the other cases, the Hecke correspondences on the modular curves
give rise to Hecke operators on the Jacobian.

\begin{prop}\label{parjac}
Let $N\ge 3$ be an integer, and $p$ a prime. Then there are 
isomorphisms of $\FF_p$-vector spaces
$$ \Hpar^1 (\Gamma_1(Np),\FF_p) \cong J(\CC)[p] = J(\Qpbar)[p]$$
respecting the Hecke operators.
\end{prop}

\pf
The second equality follows from the fact that torsion points are algebraic.
We start with the exact {\em Kummer sequence} of analytic sheaves on $X_1(Np)$
$$ 0 \to \mu_p \to \GG_m \xrightarrow{p} \GG_m \to 0.$$
Its long exact sequence in analytic cohomology yields
$$ 0 \to H^1(X_1(Np), \mu_p) \to H^1(X_1(Np),\GG_m)
\xrightarrow{p} H^1(X_1(Np),\GG_m).$$
Using $H^1(X_1(Np),\GG_m) \cong \Pic_{X_1(Np)}(\CC)$
one obtains that $H^1(X_1(Np), \mu_p)$ is isomorphic to $J(\CC)[p]$.
As $\CC$ contains the $p$-th roots of unity, we may replace
the sheaf $\mu_p$ by the constant sheaf~$\FF_p$.
Moreover, the group $H^1(X_1(Np), \FF_p)$ coincides with
$\Hpar^1(\Gamma_1(Np), \FF_p)$. This follows for example
from the Leray spectral sequence for the open immersion
$Y_1(Np) \hookrightarrow X_1(Np)$ (see e.g.\ \cite{ArtEins}).
Since the Hecke operators come from correspondences on modular curves,
the isomorphisms are compatible for the Hecke action.
\qed

\section {Hecke algebras}

In this section we compare the Hecke algebra of cuspidal modular forms
to that of parabolic group cohomology and establish isomorphisms in certain
cases. The principal result is Theorem~\ref{thmdiag}, from which
the application to ordinary cusp forms (Corollary~\ref{corap})
mentioned in the introduction follows.

Whenever for a ring~$R$ we have an $R$-module~$M$, on which Hecke operators~$T_n$
act for all~$n$, we let
$$ \TT_R(M) := R[T_n \; | \; n \in \NN] \subseteq \End_R(M),$$
i.e.\ the $R$-subalgebra of the endomorphism algebra generated by
the Hecke operators.

We claim that in the situation of Hecke operators on modular forms
or group cohomology for $\Gamma_1(N)$ of weight~$k$ the diamond
operators are in $\TT_R(M)$. For $a$ coprime to~$N$, 
let $l_1, l_2$ be distinct primes which are congruent to $a$ modulo~$N$.
Then $l_i^{k-1}\langle a \rangle = T^2_{l_i} - T_{l_i^2}$, so
that with $1= rl_1^{k-1} + s l_2^{k-1}$ the claim follows.

\subsection* {The Hecke algebra of modular forms and Eichler-Shimura}

We start by stating the Eichler-Shimura theorem 
(see \cite{DiamondIm}, Theorem~12.2.2 and Proposition~12.4.10).

\begin{thm}[Eichler-Shimura]\label{thmes} 
For $k \ge 2$ and $\Gamma \le \SL_2(\ZZ)$ a congruence subgroup
there is an isomorphism of $\TT_\ZZ(S_k(\Gamma,\CC))$-modules,
the {\em Eichler-Shimura isomorphism},
$$ \Hpar^1 (\Gamma,V_{k-2}(\CC)) 
\cong S_k(\Gamma,\CC) \oplus  \overline{S_k(\Gamma,\CC)}.$$
\end{thm}

\begin{cor}\label{cores}
For $k \ge 2$ and $N \ge 3$ we have natural ring isomorphisms 
$$ \TT_\ZZ\big(S_k(\Gamma_1(N),\CC)\big) 
\cong \TT_\ZZ\big(\Hpar^1(\Gamma_1(N),V_{k-2}(\ZZ))/\torsion\big).$$
\end{cor}

\pf
The free $\ZZ$-module 
$\Hpar^1 (\Gamma_1(N),V_{k-2}(\ZZ))/\torsion$
is a $\ZZ$-structure in the $\CC$-vector space $\Hpar^1 (\Gamma_1(N),V_{k-2}(\CC))$.
Any $\ZZ$-structure gives an isomorphic Hecke algebra.
Finally, Theorem~\ref{thmes} implies that
the Hecke algebra of $\Hpar^1 (\Gamma_1(N),V_{k-2}(\CC))$
is isomorphic to the Hecke algebra of $S_k(\Gamma_1(N),\CC)$.
\qed
\smallskip

The formula in this corollary is the reason why many people prefer
to factor out the torsion of modular symbols.

\begin{prop}\label{propheckemod}
Let $N \ge 5$, $k \ge 2$ integers and $p \nmid N$ a prime.
Then we have
$$ \TT_\ZZ \big(S_k(\Gamma_1(N),\CC)\big) \otimes_\ZZ \FF_p \cong
   \TT_{\FF_p} \big(S_k(\Gamma_1(N),\FF_p)\big).$$
\end{prop}

\pf
By \cite{DiamondIm}, Theorem 12.3.2, 
there is no difference between Katz cusp forms over~$\FF_p$
and those that are reductions of classical cusp forms 
whose $q$-expansion is in~$\ZZ[1/N]$, i.e.\
$$ S_k(\Gamma_1(N),\ZZ[1/N]) \otimes_{\ZZ[1/N]}  \FF_p \cong
   S_k(\Gamma_1(N),\FF_p).$$
Hence, the $q$-expansion principle gives the two perfect pairings
$$ \TT_{\ZZ} \big(S_k(\Gamma_1(N),\CC)\big) \otimes_\ZZ \ZZ[1/N] \times
S_k(\Gamma_1(N),\ZZ[1/N]) \to \ZZ[1/N], \;\; (T,f) \mapsto a_1(Tf)$$
and
$$ \TT_{\FF_p} \big(S_k(\Gamma_1(N),\FF_p)\big)  \times
S_k(\Gamma_1(N),\FF_p) \to \FF_p, \;\; (T,f) \mapsto a_1(Tf).$$
Tensoring the first one with $\FF_p$ allows us to compare it to the 
second one, from which the proposition follows.
\qed

\begin{cor}\label{coreseins}
Let $p$ be a prime and $N \ge 5$, $2 \le k \le p+2$ integers s.t.~$p \nmid N$.
Then sending the operator $T_l$ to $T_l$ for all primes~$l$
defines a surjective $\FF_p$-algebra homomorphism
$$ \TT_{\FF_p}\big(S_k(\Gamma_1(N),\FF_p)\big) \twoheadrightarrow
   \TT_{\FF_p}\big(\Hpar^1(\Gamma_1(N),V_{k-2}(\FF_p))\big).$$
\end{cor}

\pf
From Corollary~\ref{cores} we obtain because of $p$-torsion-freeness
(Proposition~\ref{torsion})
an isomorphism of $\FF_p$-algebras
$$\TT_\ZZ\big(S_k(\Gamma_1(N),\CC)\big) \otimes \FF_p
\cong \TT_{\ZZ}\big(\Hpar^1(\Gamma_1(N),V_{k-2}(\ZZ))\big) \otimes_{\ZZ} \FF_p.$$
By Proposition~\ref{propheckemod} the term on the left hand side is
equal to $\TT_{\FF_p}\big(S_k(\Gamma_1(N),\FF_p)\big)$ so that it suffices
to have a surjection
$$ \TT_{\ZZ}\big(\Hpar^1(\Gamma_1(N),V_{k-2}(\ZZ))\big) \otimes \FF_p
\twoheadrightarrow \TT_{\FF_p}\big(\Hpar^1(\Gamma_1(N),V_{k-2}(\FF_p))\big),$$
which follows from Proposition~\ref{torsion}. Indeed, the isomorphism
$$\Hpar^1(\Gamma_1(N),V_{k-2}(\ZZ)) \otimes \FF_p 
\cong \Hpar^1(\Gamma_1(N),V_{k-2}(\FF_p))$$
is compatible with Hecke operators, and allows one to define a homomorphism
from the Hecke algebra on the left hand term to the one on the 
right hand term, which is automatically surjective by the definition
of the Hecke algebra.
\qed

\begin{prop}\label{tepsilon}
Let $N \ge 1$, $k \ge 2$ be integers and $K$ a field.
If the characteristic of~$K$ is $p >0$, then we assume $p \nmid N$.
Furthermore, let $\epsilon: \Gamma_0(N) \proj \Gamma_0(N) / \Gamma_1(N) \to K^*$ 
be a character such that $\epsilon(-1)=(-1)^k$.
Denote by $\TT$ the $K$-Hecke algebra of $S_k(\Gamma_1(N),K)$ and by
$\TT_\epsilon$ the $K$-Hecke algebra of $S_k(\Gamma_1(N),\epsilon,K)$.
Furthermore, let
$$ I = ( \langle \delta \rangle - \epsilon (\delta) \; | \; 
\delta \in \Gamma_0(N) / \Gamma_1(N)) \lhd \TT.$$
Then $\TT/I$ and $\TT_\epsilon$ are isomorphic $K$-algebras.
\end{prop}

\pf
As we work with Katz modular cusp forms (for that we need the condition $p \nmid N$), 
we dispose of the $q$-expansion principle. Hence we have
isomorphisms respecting the Hecke action
$(\TT_\epsilon)^\vee \cong S_k(\Gamma_1(N),\epsilon,K) \cong \TT^\vee [I] \cong (\TT/I)^\vee,$
whence the proposition follows.
\qed

\subsection* {Comparing Hecke algebras over~$\FF_p$}

\begin{prop}\label{keinskzwei}
Let $N\ge 5$ be an integer, $p \nmid N$ a prime and $0 \le d \le p-1$ an integer.
There exists a surjection 
$\TT_{\FF_p}\big(\Hpar^1(\Gamma_1(Np),\FF_p)(d)\big) \twoheadrightarrow \TT_{\FF_p} (\Lbar(d))$
such that the diagram of $\FF_p$-algebras
$$ \xymatrix@=.3cm{
& \TT_{\FF_p} (\Lbar(d)) \\
\TT_\ZZ \big(S_2(\Gamma_1(Np),\CC)(d)\big) \otimes \FF_p \ar@{->>}[ru] \ar@{->>}[rd]  & \\
& \TT_{\FF_p}\big(\Hpar^1(\Gamma_1(Np),\FF_p)(d)\big) \ar@{->>}[uu] }$$
commutes. All maps are uniquely determined by sending the Hecke operator $T_l$ to $T_l$
for all primes~l.
\end{prop}

\pf
Let us first remark how the diagonal arrows are made. The lower one comes
from the isomorphism (see Proposition~\ref{torsion})
$$\Hpar^1(\Gamma_1(Np),\ZZ) \otimes \FF_p \cong \Hpar^1(\Gamma_1(Np),\FF_p).$$
The upper one is due to the fact that $L$ is a lattice in 
$S_2(\Gamma_1(Np),\QQ_p(\zeta_p))$ (see \cite{Gross}, p. 472) and
using arguments as in Corollary~\ref{cores}.
As the order of $\FF_p^*$ is invertible in~$\FF_p$, we can everywhere
pass to the eigencomponents of the action of the $p$-part of the diamond 
operator $\diap \cdot$.

We obtain the vertical arrow by showing that the kernel of the lower diagonal
map is contained in the kernel of the upper diagonal map. In other words,
we will show that if a Hecke operator $T$ in 
$\TT_\ZZ \big(S_2(\Gamma_1(Np),\CC)(d)\big) \otimes \FF_p$
acts as zero on $\Hpar^1(\Gamma_1(Np),\FF_p)(d)$, then it acts as zero on
$\Lbar(d)$.

So assume that $T$ acts as zero on $\Hpar^1(\Gamma_1(Np),\FF_p)(d)$.
By Proposition~\ref{parjac}, it acts as zero on $J_{\Qpbar}(\Qpbar)[p](d)$,
hence on $J_{\Qpbar}[p](d)$, as $J_{\QQ_p}[p]$ is reduced.
But then it also acts as zero on $J_{\ZZ_p[\zeta_p]}[p](d)$, as it acts as
zero on the generic fibre
using that $J[p]$ is flat over $\ZZ_p[\zeta_p]$ (\cite{BLR}, Lemma~7.3.2,
as $J$ is semi-abelian).
But consequently, it also acts as zero on 
the special fibre $J_{\FF_p}[p](d)$, whence also on the cotangent space
$\Cot_0(J_{\FF_p}^0[p])(d)$. Now Lemma~\ref{cotlem} finishes the proof.
\qed

\begin{thm}\label{thmdiag}
Let $p$ be a prime and $2 < k \le p+1$, $N \ge 5$ integers such that $p \nmid N$.
We write for short
$\TT^{\mathrm{par},N,k} := \TT_{\FF_p} \big(\Hpar^1(\Gamma_1(N),V_{k-2}(\FF_p))\big)$,
$\TT^{\mathrm{mod},N,k} := \TT_{\FF_p} \big(S_k(\Gamma_1(N),\FF_p)\big)$ and
similarly for the twisted ones.
Then there is the commutative diagram of $\FF_p$-algebras
$$ \xymatrix@=.7cm{
  \TT_{\FF_p} \big(\Lbar(k-2)\big) \ar@{->}[r] \ar@{<-}[d] 
& \TT^{\mathrm{mod},N,k} \ar@{->}[d] \ar@{}|(.35){\times}[r]
& \TT^{\mathrm{mod},N,p+3-k,[k-2]} \ar@{->}[d] \\
  \TT^{\mathrm{par},Np,2} \ar@{->}[r] 
& \TT^{\mathrm{par},N,k}  \ar@{}|(.35){\times}[r]
& \TT^{\mathrm{par},N,p+3-k,[k-2]}.}$$
The vertical arrows are obtained from Proposition~\ref{keinskzwei}
resp.\ Corollary~\ref{coreseins}, and the horizontal
ones from Proposition~\ref{propgross} and Proposition~\ref{heckepar}.
The vertical arrows are surjective.
If $2 < k \le p$, then the upper horizontal arrow is injective.
\end{thm}

\pf
The commutativity is clear, as $T_l$ is sent to $T_l \times T_l$ along the horizontal
arrows, and $T_l$ is sent to $T_l$ along the vertical arrows for all primes~$l$.
The surjectivity of the vertical arrows has been proved at the places
cited above.
The injectivity of the upper homomorphism is the fact that $\Lbar(k-2)$
is the direct sum of $S_k(\Gamma_1(N),\FF_p)$ and $S_{p+3-k}(\Gamma_1(N),\FF_p)[k-2]$,
if $2 < k \le p$.
\qed

\begin{cor}\label{corsoev}
Let $2 < k \le p+1$, $N \ge 5$ such that $p \nmid N$.
Let $\fP$ be a maximal ideal of the Hecke algebra
$\TT_\ZZ \big(S_2(\Gamma_1(Np),\CC)(d)\big) \otimes \FF_p$
which is not in the support of $S_{p+3-k}(\Gamma_1(N),\FF_p)$.
Then we have an isomorphism
$$ \TT_{\FF_p}\big(S_k(\Gamma_1(N),\FF_p)_\fP\big) \cong 
   \TT_{\FF_p}\big(\Hpar^1(\Gamma_1(N),V_{k-2}(\FF_p))_\fP \big).$$
\end{cor}

\pf
The assumption means that $(S_{p+3-k}(\Gamma_1(N),\FF_p)[k-2])_\fP = 0$.
By Corollary~\ref{coreseins} the ideal $\fP$ is not in the
support of $\Hpar^1(\Gamma_1(N),V_{p+1-k}(\FF_p))[k-2]$ either, whence
$(\Hpar^1(\Gamma_1(N),V_{p+1-k}(\FF_p))[k-2])_\fP = 0$.
Hence, the sequences of Propositions~\ref{heckepar} and~\ref{propgross}
localised at~$\fP$ give isomorphisms
$\TT_{\FF_p} \big(\Lbar(k-2)\big)_\fP \cong \TT^{\mathrm{mod},N,k}_\fP$
and
$\TT^{\mathrm{par},Np,2}_\fP \cong \TT^{\mathrm{par},N,k}_\fP$.
Hence, also the vertical maps in the localisation at~$\fP$ of the diagram 
of Theorem~\ref{thmdiag} are isomorphisms.
\qed

\begin{cor}\label{corap}
Let $2 < k \le p+1$, $N \ge 5$ such that $p \nmid N$.
Let $\fP$ be a maximal ideal of $\TT_{\FF_p}\big(S_k(\Gamma_1(N),\FF_p)\big)$
corresponding to a normalised eigenform $f \in S_k(\Gamma_1(N),\Fpbar)$
which is {\em ordinary}, i.e.\ $a_p(f) \neq 0$.
Then we have an isomorphism
$$ \TT_{\FF_p}\big(S_k(\Gamma_1(N),\FF_p)_\fP\big) \cong 
   \TT_{\FF_p}\big(\Hpar^1(\Gamma_1(N),V_{k-2}(\FF_p))_\fP \big).$$
\end{cor}

\pf
As the operator $T_p$ always acts as zero on $S_{p+3-k}(\Gamma_1(N),\FF_p)[k-2]$
the maximal ideal~$\fP$ cannot be in the support of $S_{p+3-k}(\Gamma_1(N),\FF_p)[k-2]$,
whence we are in the situation of Corollary~\ref{corsoev}.
\qed

\section {Application to $\Gamma_0(N)$ and characters}

The techniques from the preceding sections do not apply to the group
$\Gamma_0(N)$. In this section we show that one can, nevertheless, 
derive similar results for that group together with a 
character of the quotient group $\Gamma_0(N) / \Gamma_1(N) \cong (\ZZ/N\ZZ)^*$.

Let us for the sequel of this section make the following assumption.

\begin{ass}\label{annahme}
Let $p \ge 5$ be a prime and $K$ a finite field of characteristic~$p$
or $K = \Fpbar$. Suppose, moreover, that we are given integers $N \ge 5$ and $k$
with $3 \le k \le p+1$ and $p \nmid N$. 
Let $\Delta = \Gamma_0(N) / \Gamma_1(N)$.
It acts on $H^1(\Gamma_1(N),\cdot)$, the parabolic subspace and on 
$S_k(\Gamma_1(N),K)$ through the diamond operators.
Furthermore, let $\epsilon$ be a character of the form
$\epsilon: \Gamma_0(N) \proj \Gamma_0(N)/\Gamma_1(N) \to K^*$. 
Denote by $K^\epsilon$ the $K[\Gamma_0(N)]$-module which is a copy of~$K$ with
action through $\epsilon^{-1}$.
We write $V_{k-2}^\epsilon$ for the $K[\Gamma_0(N)]$-module 
$V_{k-2}(K) \otimes_K K^\epsilon$.
Finally, let $G$ be the group with $\Gamma_1(N) \le G \le \Gamma_0(N)$
such that $G/\Gamma_1(N) = \Delta_p$, the $p$-Sylow subgroup of~$\Delta$.
Note that $G$ acts freely on~$\HH$.
\end{ass}

\begin{lem}\label{gpdelta}
Let us assume~\ref{annahme}. Then 
$H^1(\Gamma_1(N), V_{k-2}(K))$ and its parabolic subspace
are coinduced $K[\Delta_p]$-modules.
\end{lem}

\pf
We write $V = V_{k-2}(K)$ and $\Gamma := \Gamma_1(N)$.
The Hochschild-Serre spectral sequence (in the notation from 
\cite{Milne}, Appendix~B,
we use the exact sequence $E_1^2 \to E_2^{1,1} \to E_2^{3,0}$)
gives the exact sequence
\begin{multline*}
 \ker\big( H^2(G,V) \to H^0(\Delta_p,H^2(\Gamma,V)) \big)
\\ \to H^1(\Delta_p,H^1(\Gamma,V)) \to H^3(\Delta_p,H^0(\Gamma,V)),
\end{multline*}
whence $H^1(\Delta_p,H^1(\Gamma, V))$ is zero by 
Corollary~\ref{corhzweinull} and Proposition~\ref{vknull}.
From Proposition~\ref{leraypar} it follows taking $\Delta_p$-cohomology
that the sequence
\begin{align*}
0 &\to \Hpar^1(\Gamma, V)^{\Delta_p} \to H^1(\Gamma,V)^{\Delta_p} \xrightarrow{\res} 
\big(\prod_{g \in \Gamma \backslash \PSL_2(\ZZ) / \langle T \rangle}
H^1(\Gamma \cap \langle g T g^{-1} \rangle, V)\big)^{\Delta_p} \\
& \to H^1(\Delta_p,\Hpar^1(\Gamma,V)) \to 0 
\end{align*}
is exact. It can be checked (using that $-1 \not \in G$ and 
that $p$ does not divide the ramification
indices of the cusps) that the first three terms are
$$ 0 \to \Hpar^1(G, V) \to H^1(G,V) \xrightarrow{\res} 
\prod_{g \in G \backslash \PSL_2(\ZZ) / \langle T \rangle}
H^1(G \cap \langle g T g^{-1} \rangle, V) \to 0, $$
where the zero on the right is again a consequence of Proposition~\ref{leraypar}.
Hence, $H^1(\Delta_p,\Hpar^1(\Gamma, V))$ is zero.
Finally, \cite{NSW}, Proposition 1.7.3(ii), implies that
$H^1(\Gamma,V)$ and its parabolic subspace
are coinduced $\Delta_p$-modules.
\qed

\begin{lem}\label{moddelta}
Let us assume~\ref{annahme}. Then $S_k(\Gamma_1(N),K)$
is a free $K[\Delta_p]$-module.
\end{lem}

\pf
For the notation in this proof we follow \cite{EdixBoston}, pp.~209-210
and the proof of Lemma~1.9. We let $\Gamma := \Gamma_1(N)$.

The projection 
$\pi: X_\Gamma \twoheadrightarrow X_G$ is a Galois 
cover with group~$\Delta_p$ of projective $K$-schemes. Indeed, for the
open part this is \cite{DR}, VI.2.7. Moreover,
the ramification index of the cusps divides~$N$,
whence the cusps are unramified in a $p$-extension.

Using Serre duality and the Kodaira-Spencer isomorphism 
(see \cite{DR}, VI.4.5.2) we obtain
\begin{align*}
 H^1 (X_G, \omega^{\otimes k}(-\cusps)) & \overset{\text{S-D}}{\cong}
H^0 (X_G, \Omega^1 \otimes ( \omega^{\otimes k} (-\cusps))^\vee)^\vee \\
& \overset{\text{K-S}}{\cong}
H^0(X_G, \omega^{\otimes 2-k})^\vee
\end{align*}
which is zero, since the degree of $\omega^{\otimes 2-k}$ is negative (as $k \ge 3$).
The map $\pi$ is \'etale and we have 
$H^0(X_\Gamma, \pi^* \omega^{\otimes k}(-\cusps)) \cong S_k(\Gamma,K)$.
We conclude from \cite{Nakajima}, Theorem~2, that this is a
free $K[\Delta_p]$-module.
\qed

\begin{thm}\label{heckedelta}
Let us assume \ref{annahme} and that 
$\Hpar^1(\Gamma_1(N), V_{k-2}(K))$ is a faithful $\TT_K(S_k(\Gamma_1(N),K))$-module.
Then $\Hpar^1(\Gamma_0(N),V_{k-2}^\epsilon)$ is a faithful 
module for $\TT_K(S_k(\Gamma_1(N),\epsilon,K))$.
\end{thm}

\pf
Let $\Gamma := \Gamma_1(N)$.
We claim that $N_\Delta := \sum_{\delta \in \Delta} \delta \in K[\Delta]$ 
induces isomorphisms
\begin{equation}\label{eqsk}
(S_k(\Gamma,K) \otimes_K K^\epsilon)_\Delta \to (S_k(\Gamma,K) \otimes_K K^\epsilon)^\Delta
\end{equation}
and
\begin{equation}\label{eqvk}
\Hpar^1(\Gamma,V_{k-2}^\epsilon)_\Delta \to \Hpar^1(\Gamma,V_{k-2}^\epsilon)^\Delta.
\end{equation}
We note that $H^1(\Gamma,V_{k-2}(K)) \otimes_K K^\epsilon = H^1(\Gamma,V_{k-2}^\epsilon)$, 
since the character~$\epsilon$ restricted to~$\Gamma$ is trivial.
Using Lemmas \ref{gpdelta} and \ref{moddelta},
an elementary calculation gives the claim.

Dualising Equation~(\ref{eqsk}) gives an isomorphism
$$\big( \TT(S_k(\Gamma, K)) \otimes K^\epsilon \big)_\Delta 
  \xrightarrow{N_\Delta} \big( \TT(S_k(\Gamma, K)) \otimes K^\epsilon \big)^\Delta,$$
which in particular yields the implication
\begin{equation}\label{eqtn}
T (\sum_{\delta \in \Delta} \epsilon(\delta)^{-1} \langle \delta \rangle) = 0 \;\;\;
\Rightarrow \;\;\; T \in I,
\end{equation}
where $I$ is the ideal defined in Proposition~\ref{tepsilon}. In view of that
proposition, we only need to show that if $T$ acts as zero on $H^1(G,V_{k-2}^\epsilon)$,
then $T$ is in~$I$.

The Hochschild-Serre spectral sequence yields 
$H^1(\Gamma,V_{k-2}^\epsilon)^\Delta \cong H^1(G,V_{k-2}^\epsilon)$,
since $(V_{k-2}^\epsilon)^\Gamma$ is zero by Proposition~\ref{vknull}.
Let now $T$ be a Hecke operator. Then we have
\begin{align*}
T \cdot H^1(G,V_{k-2}^\epsilon) & = T \cdot H^1(\Gamma,V_{k-2}^\epsilon)^\Delta \\
& = T N_\Delta \cdot H^1(\Gamma,V_{k-2}^\epsilon)_\Delta = 
T (\sum_{\delta \in \Delta} \epsilon(\delta)^{-1} \langle \delta \rangle) \cdot H^1(\Gamma,V_{k-2}(K))
\end{align*}
Suppose that this is zero. Hence, by the assumed faithfulness we have that
$T (\sum_{\delta \in \Delta} \epsilon(\delta)^{-1} \langle \delta \rangle) = 0$, 
which by Equation~(\ref{eqtn}) implies $T \in I$, as required.
\qed

\begin{rem}
If $k=2$, then the statements of Theorem~\ref{heckedelta} also hold
``outside the Eisenstein part''. The Eisenstein part is the subspace
on which the Hecke algebra acts via a system of eigenvalues that does not
belong to an irreducible Galois representation.
For, a suitable analogue of Lemma~\ref{moddelta} holds, since
$H^1 (X_G, \omega^{\otimes k}(-\cusps)) \cong H^0(X_G,\cO)^\vee$
cannot give rise to a non-Eisenstein system of eigenvalues.
Moreover, under the extra assumption
also Lemma~\ref{gpdelta} remains valid, as an easy calculation
shows.
\end{rem}

\section{Application to weight one modular forms}

Edixhoven explains in \cite{EdixJussieu}, Section~4, how weight one cuspidal 
Katz modular forms over finite fields of characteristic~$p$ can be computed from the
knowledge of the Hecke algebra of weight~$p$ cusp forms over the same field.
In this section we shall first recall this and then derive the conclusion
mentioned in the introduction that weight one modular cusp forms over~$\Fpbar$
can be computed by the Hecke operators on weight~$p$ parabolic group 
cohomology over~$\FF_p$.

Let $\FF$ be a finite field of prime characteristic~$p$ or $\Fpbar$ and 
fix a level~$N\ge 1$ with $p \nmid N$ and a character
$\epsilon: (\ZZ/N\ZZ)^* \to \FF^*$ with $\epsilon(-1)=(-1)^k$. 
We have two injections of $\FF$-vector spaces
$$ F,A: S_1 (\Gamma_1(N),\epsilon,\FF) 
   \to S_p (\Gamma_1(N),\epsilon,\FF), $$
given on $q$-expansions by $a_n(Ag) = a_n(g)$ and 
$a_n(Fg) = a_{n/p}(g)$ (with $a_n(Fg) = 0$ if $p \nmid n$),
which are compatible with all Hecke operators~$T_l$ for primes $l \neq p$. 
The former comes from the {\em Frobenius} and the latter is multiplication
by the {\em Hasse invariant}. 
One has $T_p^{(p)} F = A$ and 
$A T_p^{(1)} = T_p^{(p)} A + \epsilon(p)F$, where we have indicated the
weight as a superscript (see e.g.\ \cite{EdixJussieu}, Equation~(4.1.2)).

Let $\TT^{(k)}$ be the Hecke algebra over~$\FF$
of weight~$k$ for a fixed level~$N$ and a fixed character~$\epsilon$. 
We will also indicate the weight of Hecke operators by superscripts. 
We denote by $A^{(p)}$ the $\FF_p$-subalgebra of $\TT^{(p)}$
generated by all Hecke operators $T_n^{(p)}$ for $p \nmid n$. 

\begin{prop}\label{thetaprop}
\begin{enumerate}[(a)]
\item There is a homomorphism $\Theta$, called a {\em derivation},
which on $q$-expansions is given by
$a_n(\Theta f) = n a_n(f)$ such that the sequence
$$ 0 \to S_1 (\Gamma_1(N),\epsilon,\FF) \xrightarrow{F} 
         S_p (\Gamma_1(N),\epsilon,\FF) \xrightarrow{\Theta}
         S_{p+2} (\Gamma_1(N),\epsilon,\FF) $$
is exact. 

\item Suppose $f \in  S_1(\Gamma_1(N),\epsilon,\FF)$ such that $a_n(f) = 0$ for
all~$n$ with $p \nmid n$. Then $f = 0$. 
In particular $ A S_1(\Gamma_1(N),\epsilon,\FF) \cap  F S_1(\Gamma_1(N),\epsilon,\FF) = 0$. 

\item The Hecke algebra $\TT^{(1)}$ in weight one can be generated by
all $T_l^{(1)}$, where $l$ runs through the primes different from~$p$. 

\item The weight one Hecke algebra $\TT^{(1)}$ is the algebra generated
by the $A^{(p)}$-action on the module $\TT^{(p)}/A^{(p)}$. 

\end{enumerate}
\end{prop}

\pf
(a) The main theorem of \cite{KatzDerivation} gives the exact sequence
$$ 0 \to S_1 (\Gamma_1(N),\epsilon,\FF) \xrightarrow{F} 
         S_p (\Gamma_1(N),\epsilon,\FF) \xrightarrow{A\Theta}
         S_{2p+1} (\Gamma_1(N),\epsilon,\FF) $$
by taking Galois invariants.
However, as explained in \cite{EdixJussieu}, Section~4, 
the image $A\Theta S_p(\Gamma_1(N),\epsilon,\FF)$ in weight $2p+1$ can be
divided by the Hasse invariant, whence the weight is as claimed. 

(b) The condition implies by looking 
at $q$-expansions that $A\Theta f = 0$, whence by Part~(3) of Katz' theorem
cited above $f$ comes from a lower weight than~$1$, but below there is just the 
$0$-form (see also \cite{EdixJussieu}, Proposition 4.4). 

(c) It is enough to show that $T_p^{(1)}$ is linearly dependent on the span of all
$T_n^{(1)}$ for $p \nmid n$. If it were not, then there would be a modular cusp form
of weight~$1$ satisfying $a_n(f) = 0$ for $p \nmid n$,
but $a_p(f) \neq 0$, contradicting~(b). 

(d) Dualising the exact sequence in~(a) yields that $\TT^{(p)}/A^{(p)}$ and $\TT^{(1)}$
are isomorphic as $A^{(p)}$-modules, which implies the claim. 
\qed

\begin{prop}\label{HeckeBound}
Let $N \ge 1$ and $k \ge 2$ be integers such that $p \nmid N$, 
$\FF | \FF_p$ a finite extension
and let $\epsilon: (\ZZ/N\ZZ)^* \to \FF^*$ be a character with $\epsilon(-1)=(-1)^k$. 
Set
$$ B = \frac{N}{12}\prod_{l \mid N, l \text{ prime}} (1 + \frac{1}{l}).$$
\begin{enumerate}[(a)]
\item Then the Hecke operators
$T_1^{(k)}$, $T_2^{(k)}$, $\dots$, $T_{kB}^{(k)}$
generate $\TT^{(k)}$ as an $\FF$-vector space. 

\item The $\FF$-algebra $A^{(p)}$
can already be generated as an $\FF$-vector space by the set
$$\{\; T_n^{(p)} \; \mid \; p \nmid n, \; n \le (p+2)B \; \}.$$
\end{enumerate}
\end{prop}

\pf
(a) This follows from the proof of \cite{EdixJussieu}, Proposition~4.2. 

(b) Assume that some
$T_m^{(p)}$ for $m > (p+2)B$ and $p \nmid m$ is linearly independent of the
operators in the set of the assertion. This means that there
is a cusp form $f \in S_p(\Gamma_1(N),\epsilon,\FF)$
satisfying $a_n(f) = 0$ for all $n \le (p+2)B$ with $p \nmid n$, but $a_m(f) \neq 0$. 
One gets $a_n(\Theta f) = 0$ for all $n \le (p+2)B$, but
$a_m(\Theta f) \neq 0$. This contradicts~(a). 
\qed

\begin{rem}\label{baccent}
If we work with $\Gamma_1(N)$ and no character, the number~$B$ above
has to be replaced by
$$ B' = \frac{N^2}{24}\prod_{l \mid N, l \text{ prime}} (1 - \frac{1}{l^2}).$$
\end{rem}

Part of the following proposition is \cite{EdixJussieu}, Proposition~6.2.

\begin{prop}\label{eigenspace}
Let $V \subset S_p(\Gamma_1(N),\epsilon,\Fbar)$ be the eigenspace of a
system of eigenvalues for the operators $T_l^{(p)}$ for all primes $l \neq p$

If the system of eigenvalues does not come from a weight one form,
then $V$ is at most of dimension one.
Conversely, if there is a normalised weight one eigenform~$g$ 
with that system of eigenvalues
for $T_l^{(1)}$ for all primes $l \neq p$, then $V = \langle Ag, Fg \rangle$ and that
space is $2$-dimensional. 
On it $T_p^{(p)}$ acts with eigenvalues $u$ and $\epsilon(p)u^{-1}$
satisfying $u+\epsilon(p)u^{-1} = a_p(g)$. 
In particular, the eigenforms in weight~$p$ which come from weight one are ordinary. 
\end{prop}

\pf
We choose a normalised eigenform~$f$ for all operators.
If $V$ is at least~$2$-dimensional, then we have 
$V = \FF f \oplus \{h \;|\; a_n(h) = 0 \; \forall p \nmid n\}$. 
As a form $h$ in the right summand is annihilated by~$\Theta$, it is
equal to $Fg$ for some form~$g$ of weight one by Proposition~\ref{thetaprop}~(a). 
By Part~(b) of that proposition we know that $\langle Ag, Fg \rangle$ 
is $2$-dimensional. If $V$ were more than~$2$-dimensional, then there
would be two different cusp forms in weight~$1$, which are eigenforms for all $T_l^{(1)}$
with $l \neq p$. This, however, contradicts Part~(c). 

Assume now that $V$ is $2$-dimensional.
Any normalised eigenform $f \in V$ for all Hecke operators in weight~$p$
has to be of the form $Ag + \mu Fg$ for some~$\mu \in \Fbar$. 
The eigenvalue of $T_p^{(p)}$ on~$f$ is the $p$-th coefficient,
hence $u = a_p(g) + \mu$, as $a_p(Fg) = a_1(g) = 1$. 
Now we have
\begin{align*}
 (a_p(g) + \mu) (Ag + \mu Fg)     
  & =  T_p^{(p)} (Ag+\mu Fg)    = T_p^{(p)}Ag + \mu Ag \\
  & =  A T_p^{(1)} g - \epsilon(p)Fg + \mu A g    =  (a_p(g) + \mu) A g - \epsilon(p)Fg,
\end{align*} 
which implies $- \epsilon(p) = (a_p(g) + \mu)\mu = u^2- ua_p(g)$ 
by looking at the $p$-th coefficient. From this one obtains the claim on~$u$. 
\qed

\begin{thm}\label{wtoneOK}
Let $N \ge 5$ an integer and $p$ be a prime not dividing~$N$. 

\begin{enumerate}[(a)]
\item The Hecke algebra of
$S_1(\Gamma_1(N),\FF_p)$ can be computed using 
the first $(p+2)B'$ Hecke operators on $\Hpar^1(\Gamma_1(N),V_{p-2}(\FF_p))$.
\item Let $\FF$ be a finite field of characteristic $p \ge 5$ 
and let $\epsilon: (\ZZ/N\ZZ)^* \to \FF^*$ be a character. 
The Hecke algebra of $S_1(\Gamma_1(N),\epsilon,\FF)$ can be computed using 
the first $(p+2)B$ Hecke operators on $\Hpar^1(\Gamma_0(N),V_{p-2}^\epsilon(\FF))$.
\end{enumerate}
The numbers $B$ resp.\ $B'$ were defined in Proposition~\ref{HeckeBound} 
and Remark~\ref{baccent}.
\end{thm}

\pf
(a) Corollary~\ref{corap} implies that the ordinary part of $\Hpar^1(\Gamma_1(N),\FF_p)$
is a faithful module for the ordinary part of the Hecke algebra of weight~$p$ Katz cusp
forms over~$\FF_p$. So that part of the Hecke algebra can be computed using 
the Hecke operators $T_1, \dots, T_{pB'}$ on $\Hpar^1(\Gamma_1(N),\FF_p)$ (see Proposition~\ref{HeckeBound}(a)).
From Proposition~\ref{eigenspace} we know that the image of the weight one
forms in weight~$p$ under the Hasse invariant and Frobenius lies in the ordinary part.
Proposition~\ref{HeckeBound}(b) implies that $A^{(p)}$ can be computed by the
Hecke operators indicated there.
Now the Hecke algebra of weight one Katz cusp forms on $\Gamma_1(N)$ without 
a character can be computed as described in Proposition~\ref{thetaprop}(d).

(b) Under the extra assumption and using Theorem~\ref{heckedelta},
the same arguments also work with a character and the bound~$B$.
\qed

\vspace*{.5cm}
\noindent Gabor Wiese\\
NWF I - Mathematik\\
Universit\"at Regensburg\\
D-93040 Regensburg\\
Deutschland\\
E-mail: {\tt gabor@pratum.net}\\
Web page: {\tt http://maths.pratum.net/}

\end{document}